\documentclass[11pt,a4paper]{article}

\title{Stratified simplices and intersection homology}

\author{Jonathan Fine\relax
\thanks{203 Coldhams Lane, Cambridge, CB1 3HY, England.  
\quad E-mail: \texttt{j.fine@pmms.cam.ac.uk}\hfil\break
This version: 17 September 1998}
}
\date{23 July 1998}

\newtheorem{theorem}{Theorem}
\newtheorem{definition}[theorem]{Definition}
\newcommand\bibrule{\rule{2pc}{0.4pt}}

\def\PD{{\bf P\!}_\Delta}
\def\bfR{{\bf R}}
\def\bfC{{\bf C}}
\def\bfe{{\bf e}}
\def\bfP{{\bf P}}

\begin{document}
\maketitle

\begin{abstract}
\noindent
Intersection homology is obtained from ordinary homology by imposing
conditions on how the embedded simplices meet the strata of a space $X$.
In this way, for the middle perversity, properties such as strong
Lefschetz are preserved.  This paper defines local-global intersection
homology groups, that record global information about the singularities of
$X$.  They differ from intersection homology in that stratified rather
than ordinary simplices are used.  An example of such is $\sigma_j\times
C\sigma_i$, where $\sigma_i$ and $\sigma_j$ are ordinary simplices, and
$C$ is the coning operator.  The paper concludes with a sketch of the
relationship between local-global homology and the geometry of convex
polytopes.  This paper is a more formal exposition of part of the author's
\emph{Local-global intersection homology}, alg-geom/9709011.
\end{abstract}

\section {Introduction}

Throughout $X$ will be an irreducible complex algebraic variety of complex
dimension $n$, considered as a topological space of real dimension $m=2n$. 
Homology theory, as developed by Poincar\'e, associates to $X$ a family
$H_iX$ of ordinary homology groups, where the index $i$ runs from $0$ to
$m$.  Throughout, we will calculate homology with real coefficients.  When
$X$ is nonsingular, as is well known, these homology groups have the
following properties.

First, there is defined a natural intersection pairing
\[
    H_iX \otimes H_j X \to H_0X \cong \bfR
\]
whenever $i+j=m$, and this pairing is perfect (non-degenerate).  This is
\emph{Poincar\'e duality}.  Second, each embedding $X\subset \bfP_N$ of X
in a projective space gives rise to a hyperplane class or \emph{Lefschetz
element} $\omega \in H_{m-2}X$ with the following property.  There defined
is a natural intersection map
\[
    \omega : H_iX \to H_{i-2}X
\]
that is injective for $i\geq n$ and surjective for $i-2\leq n$. 
Equivalently, for $i+j=m$, $i\geq j$ and $2k=i-j$ the map
\[
    \omega^k : H_iX \to H_jX
\]
is an isomorphism.  This is the \emph{strong Lefschetz theorem}.  Thirdly,
as a consequence of Deligne's proof of the Weil conjectures, the
\emph{Betti numbers} $h_iX$ (the dimension of $H_iX$) can be computed by
counting points on $X$ over finite fields, and thus there are
combinatorial formulae for the Betti numbers $h_iX$.

Now suppose that the algebraic variety $X$ in singular.  In this case, as
a topological space $X$ need not be a manifold, and in general the
properties listed do not hold for the ordinary homology of $X$.  They are
true however for the middle perversity intersection homology (mpih) theory
introduced by Goresky and MacPherson.  This theory is in its definition
similar to ordinary homology, except that it imposes restrictions on how
the cycles and boundaries meet the singularities.  Poincar\'e duality was
proved by Goresky and MacPherson \cite{bib.MG-RDM.IH}.  Strong Lefschetz
was proved in \cite{bib.AB-JB-PD.FP}, again as a consequence of the proof
of the Weil conjectures.  For toric varieties, see any of
\cite{bib.JD-FL.IHNP,bib.KF.IHTV,bib.RS.GHV} for the Betti number formula.
Throughout, the middle perversity will be used.

The ordinary homology of nonsingular varieties has other properties, such
as the Riemann-Hodge inequalities and a ring structure, whose analogues
for the intersection homology of singular varieties is at present
unknown.  For the author, combinatorial formulae for the Betti numbers are
the most significant property, whereas classically it is the functorial
nature that is pre-eminent.

\section{Ordinary and intersection homology}

In this section and the next we review matters, so as to provide some
background and motivation.

As is well known, the ordinary homology of $X$ can be defined in the
following way.  The \emph{standard $i$-simplex} $\sigma_i$ is the convex
hull of the standard basis vectors $\bfe_0$, $\ldots$, $\bfe_i$ in
$\bfR^{i+1}$.  Alternatively, it is the region defined by the inequalities
$x_0\geq 0$, $\ldots$, $x_i\geq 0$ and the equation $x_0 + \dots + x_i =
1$.  An \emph{embedded $i$-simplex} $f$ is a continuous map $f:\sigma \to
X$.

Each embedded $i$-simplex $f$ has a \emph{boundary} $df$, which is a
formal sum of embedded $(i-1)$-simplices.  More exactly
\[
    df = \partial_0 f - \partial_1 f + \dots + (-1)^i \partial_i f
\]
where $\partial_j f$ is the composition of $f$ with the $j$-th
\emph{inclusion map} $\sigma_{i-1} \hookrightarrow \sigma_i$.  (Each point
in $\sigma_{i-1}$ is a sequence of $i$ numbers.  This map prolongs the
sequence by inserting a zero immediately before the entry indexed by $j$,
or after the last entry if $j=i$.)  It is
easily seen that $d^2f=0$ for any embedded simplex $f$.  Later, we will
need to know what the boundary of a product of ordinary simplices is.

To continue, an \emph{$i$-chain} $\xi$ is simply a formal sum of embedded
$i$-simplices.  A chain $\eta$ is said to be a \emph{cycle} if its boundary
$d\eta$, also as a formal sum, is zero.  Because $d\circ d$ is zero, the
boundary $d\xi$ of a chain $\xi$ will always be a cycle.  A cycle $\eta$
is called a \emph{boundary} if the equation $\eta=d\xi$ can be solved, for
some chain $\xi$.  The $i$-th \emph{ordinary homology group} $H_iX$ of $X$
is defined to be the quotient group of the $i$-cycles modulo the
$i$-cycles that are boundaries.

Intersection homology is obtained by imposing conditions on the embedded
simplices.  To do this, we must first stratify $X$.  Throughout we will
suppose that $X$ has been written as a disjoint union
\[
    X = S_0 \sqcup S_1 \sqcup \dots \sqcup S_n
\]
of strata $S_j$.  At each of its points $s_j$ each $S_j$ is to have
complex dimension $j$, and about each $s_j$ the topological structure of
$X$ is to be locally constant.  The \emph{middle perversity conditions} on
an embedded $i$-simplex $f$ are this.  For each $j<n$ the dimension of
$f^{-1}(S_j)$ must be at most $(i-1)-(n-j)$.  Put another way, there is
of course no condition on $f^{-1}(S_n)$, but as $j$ goes from $n$ to $0$,
so the largest allowed dimension drops first by two, and then by one at
each step.

An \emph{allowed} or \emph{admitted} $i$-cycle $\eta$ is as before a
formal sum of $i$-chains whose boundary $d\eta$ is zero, where in addition
each embedded simplex in $\eta$ satisfies the above conditions.  Such a
cycle is an \emph{allowed boundary} if the equation $\eta=d\xi$ can be
solved, where the embedded simplices of $\xi$ satisfy the above
conditions.  The $i$-th (middle perversity) \emph{intersection homology}
group $IH_iX$ of $X$ consists of the allowed $i$-cycles modulo the allowed
boundaries.

\section{Local intersection homology}

This section continues the review of background material, and concludes
with some examples of local-global homology groups.

As mentioned earlier, the (middle perversity) intersection homology
groups have the Poincar\'e duality, strong Lefschetz and combinatorial
formula properties.  However, when $X$ is singular there can be
significant information about the singularities that is not recorded by
these groups.  This information is recorded after a fashion by
intersection homology with other perversities, or by the various `change
of perversity' groups that can be defined.  But these groups generally
fail to satisfy strong Lefschetz, and do not have a combinatorial formulae
for the Betti numbers \cite{bib.McC.HTV}.  The local-global theory is an
attempt to define additional homology groups, whose Betti numbers are we
hope given by a combinatorial formula.

To motivate the definition, we will study the topology of $X$ about a
stratum $S_j$.  First let $s$ be a point on $S_j$.  About $s$ the space
$X$ looks like the cone on something, and that something is the
\emph{link} $L_s$ of $X$ about $s$.  Up to homeomorphism, $L_s$ is locally
constant on $S_j$.  Now suppose that $\eta$ is an $i$-cycle on $L_s$, with
$i\geq 1$.  The space $X$ is about $s$ homeomorphic to $CL_s$, and so in
$X$ the cycle $\eta$ can be `coned away to the apex'.  In other words, it
is the boundary of the cone $C\eta$ of $\eta$, or $\eta=dC\eta$.  For
ordinary homology this expresses $\eta$ as a boundary, but for
intersection homology the $(i+1)$-chain $C\eta$ is not allowed.  Put
another way, ordinary homology is locally trivial, but intersection
homology is not.

In this way one can define the local intersection homology groups of $X$
at $s$.  When $X$ is a Schubert variety they have been much studied, and
their Betti numbers are given by the celebrated Kazhdan-Lusztig
polynomials.  These groups do not however record global information about
$X$.  To take a step in this direction, we will allow the point $s$ to
move.  To fix ideas, suppose that $s_t$, $0\leq t \leq 1$ is a one
dimensional path $\gamma$ on $X$.  So long as $s_t$ stays on $S_j$, local
constancy allows a local cycle $\eta_0$ at $s_0$ to be moved to another
local cycle $\eta_1$ at $s_1$.  Now suppose that $s_0$ and $s_1$ are the
same point.  Although $\eta_0$ and $\eta_1$ are equivalent along $\gamma$,
they need not be equivalent in a neighborhood of $s_0$.  This is an
example of the \emph{mondromy action of the fundamental group of $S_j$} on
the local homology at $s\in S_j$.

Now suppose that $s_t$, $0\leq t \leq 1$ is again a path $\gamma$ on $X$,
but this time we will allow $s_t$ to change strata.  More exactly, we
shall suppose that for $t > 0$ one has $s_t \in S_j$, but that $s_0$ lies
on some other strata $S_k$.  (Necessarily, $k$ will be smaller than $j$.) 
Now let $\eta_1$ be a local cycle at $s_1$.  We can again move $\eta_1$
along $\gamma$ to give local cycles $\eta_t$ at $s_t$, provided $t$ is not
zero.  When $t$ is close to zero, $\eta_t$ will lie in a small ball about
$s_0$, and so can be thought of as a cycle local to $s_0$.  Thus, paths
from $s_j$ to $s_k$ will take local cycles from $S_j$ to $S_k$.

Now form the group consisting of all formal sums of local $i$-cycles on
$X$, modulo the relations due to (a)~local boundaries, (b)~equivalence
along paths lying entirely on strata, and (c)~equivalence due to paths as
in the previous paragraph.  This is a gluing together of local homology
groups, and is a global invariant of $X$.  We will show that it does not
depend on the choice of a stratification.

First, choose a point $s_{i,j}$ on each connected component of each
stratum $S_i$ on $X$.  We are `gluing together' the local homology groups
at the $s_{i,j}$.  Now refine the stratification.  (If need be, first move
the $s_{i,j}$ so that they avoid the new components.)  This will add new
$s_{i,j}$, and may potentially reduce the relations that are available. 
But each new $s_{i,j}$ will lie on an old stratum connected component, and
so the local homology at such a new point will be equivalent to the local
homology on the old component.

Now for the relations.  Suppose that $\gamma$ is a path, as previously
considered.  Just as removing the origin does not disconnect $\bfC$, so
removing subvarieties from $X$ does not destroy the equivalence induced by
$\gamma$, provided its end points are left intact. This proves the result.

The above definition looks at equivalence of (formal sums of) local cycles
over the whole of $X$.  If instead of $X$ one applies the definition to
$S_k\cup \dots \cup S_n$, one obtains additional invariants of $X$.  The
same arguments as before demonstrates that these groups are stratification
independence.

These groups are the simplest significant examples of local-global
intersection homology.  In the next section we introduce stratified
simplices, which are the basic building blocks for local-global cycles, in
their topological form.  The following section will then define the
local-global groups, first without regard to the avoidance of strata, and
then with such regard.

\section{Stratified simplices}

Just as ordinary intersection homology can be defined using embeddings of
ordinary simplices $\sigma_i$, so local-global homology can be defined
using stratified simplices.  These objects can also be used to give a
formal expression to some of the concepts of the previous section.

For example, suppose that $\eta_1$ is an (intersection homology) $i$-cycle
local to the point $s$ of $X$, with $i\geq 1$.  By identifying a ball
about $s$ with the cone $CL_s$ on the link $L_s$ at $s$, one can produce
from $\eta_1$ a family of $\eta_t$ of local cycles, for $1\geq t> 0$, such
that in the limit $t=0$ the `cycle' $\eta_0$ is supported by $s$.  This
is, as mentioned earlier, an example of the `coning away' of a cycle. 
Such are allowed by ordinary homology, but are in general inadmissable for
intersection homology.

We will now describe this example, without explicit reference to the base
point $s$.  First let $f_1:\sigma_i\to X$ be one of the embedded simplices
out of which the cycle $\eta_1$ is constructed.  Because $\eta$ is local
(to $s$) we can extend $f_1$ to the cone $C\sigma_i$ on the simplex
$\sigma_i$, to obtain a continuous map $f:C\sigma_i \to X$.  This is an
example of an embedded stratified simplex.  Its boundary will be defined
so that the formal sum corresponding to $\eta_1$ will again be a cycle. In
other words, the boundary $df$ will be a formal sum of maps from
$C\sigma_{i-1}$.  Put another way, $f$ is a continuous function
$f(t,\lambda)$ of $t$ and $\lambda$, where $\lambda \in \sigma_i$ and
$t\in [0,1]$, subject to the constraint that $f(0,\lambda)$ is independent
of $\lambda$.  We take the boundary in the $\lambda$ directions, but not
in the $t$ direction.

Now suppose that $\eta$ is a formal sum of embedded $C\sigma_i$, such that
for each $1 \geq t >0$ the resulting $\eta_t$ is an intersection homology
$i$-cycle.  This is equivalent to saying that $d\eta$ as just described is
zero, and also that each $f(t,\lambda)$ in $\eta$ is admissable as an
embedded $i$-simplex, for $0<t\leq 1$.  Because $C\sigma_i$ is a cone it
follows that for each embedding of $C\sigma_i$ the point $s=f(0,\lambda)$
will be independent of the choice of $\lambda$ in $\sigma_i$.  In this way
the existence of `base points' is, as promised, implicit in the
construction.  It should now be clear that if $\eta_1$ is a formal sum of
local cycles, as used in the previous section, then it can be presented as
the $t=1$ surface of a formal sum $\eta$ of embedded $C\sigma_i$, and
vice versa.

The relations used between these $\eta$ are more subtle.  If $\eta_1$ and
$\eta'_1$ are equivalent as local cycles then there will, by definition,
be a formal sum of embedded $\sigma_{i+1}$, whose boundary is $\eta_1 -
\eta'_1$.  This can then be coned to give a formal sum of $C\sigma_{i+1}$,
whose boundary is $\eta-\eta'$.  Although such relations are perfectly
valid, they do not allow the base points to move.  (In the previous
section locally constancy was used to do this.)

To move $\eta_1$ at $s$ to an equivalent $\eta'_1$ at $s'$, we need a path
$\gamma$ from $s$ to $s'$.  This path is a map from $[0,1]$ to $X$.  Thus,
to produce a chain whose boundary is $\eta-\eta'$ it is natural to
consider embeddings of $\sigma_1\times C\sigma_i$, where $\sigma_1$ is of
course a $1$-simplex.  If $f$ is an embedding of $\sigma_1\times
C\sigma_i$ then its boundary will be a formal sum of embeddings of
$\sigma_0\times C\sigma_i \cong C\sigma_i$ and also of $\sigma_1\times
C\sigma_{i-1}$.  As before, we say that $f$ is \emph{allowed} if for each
$0<t\leq 1$ the embedding $f_t$ of $\sigma_1\times\sigma_i$ is allowed by
the perversity conditions.  We will now say that a formal sum $\eta$ of
embedded $C\sigma_i$ is a \emph{boundary} if first it is a chain, and
second the equation $\eta=d\xi$ can be solved, where $\xi$ is a formal sum
of allowed embedded $C\sigma_{i+1}$ and $\sigma_1\times C\sigma_i$
simplices.  That the stratified simplices used to construct $\xi$ are not
all of the same type is a subtlety not previously present.

The basic operations in the above examples are first applying the cone
operator $C$ to a simplex, and second multiplying by an ordinary
simplex.  We can now give the main definition of this section.

\begin{definition}
An \emph{order zero stratified simplex} is just an ordinary simplex
$\sigma$.  If $\sigma$ is an order $(r-1)$ stratified simplex and
$\sigma_i$ is an ordinary simplex (of dimension $i$) then $\sigma_i\times
C\sigma$ is an \emph{order $r$ stratified simplex}, and all such arise in
this way, for $r\geq 1$.

The \emph{dimension sequence} (or \emph{dimension} for short) of an
ordinary $i$-simplex is $(i)$.  The dimension of $\sigma_i\times C\sigma$
is $(\ldots, i)$, where the dots denote the entries of the dimension
sequence of $\sigma$.  We will use $\sigma(i_0,\ldots, i_r)$ (or more
briefly $\sigma(i)$ or even just $\sigma_i$) to denote the standard
stratified simplex of dimension $i=(i_0,\ldots,i_r)$.
\end{definition}

\section{Local-global intersection homology}

The usual concepts, when applied to stratified simplices, will produce
some but not all of the local-global intersection homology groups.  The
remainder arise through the imposition of additional conditions, as to
where the `apex locii' of the embedded simplices meet the strata of $X$.

Throughout this section $\sigma$ will be an order $r$ stratified simplex
of dimension $i=(i_1,\ldots,i_r)$.  As usual, an \emph{embedding} of
$\sigma$ is just a continuous map $f:\sigma\to X$.  We will say that is is
\emph{allowed} if, whenever the coning variables $t$ are all nonzero, the
restriction
\[
    f_t : \sigma (i_r) \times \dots \times \sigma(i_0) \to X
\]
of $f$ to the $t$-slice is allowed by the perversity conditions.

The \emph{boundary} $df$ of $f$ is defined so that the $t$-slice $(df)_t$
of $df$ is equal to the boundary $d(f_t)$ of the $t$-slice of $f$.  This
condition is sufficient to determine $df$.  It also ensures that $d \circ
d$ is zero.  Put another way, the construction is this.  The product of
simplices $\sigma(i_r)\times \dots \times \sigma(i_0)$ has $i_0 + \dots +
i_r + r$ boundary facets, each one of which when suitably coned up
produces a facet of $\sigma$.  (This assumes that none of the $i_j$
are zero.  But such can be ignored.)  Each facet of a $\sigma(i_j)$
determines a facet of $\sigma$.  The sign associated to such a facet is to
be determined by its location (even or odd) in $\sigma(i_r)\times \dots
\times\sigma(i_0)$, and not by its location in $\sigma(i_j)$ itself.  This
ensures in the usual way that $d\circ d$ is zero.  The \emph{boundary}
$df$ of $f:\sigma\to X$ is the signed formal sum of the restriction of $f$
to these boundary facets.

As usual, an \emph{$i$-cycle} is a formal sum $\eta$ of allowed embedded
$i$-simplices, whose boundary $d\eta$ is zero.  Such a cycle is a
\emph{boundary} if the equation $\eta=d\xi$ has a solution, where $\xi$ is
a formal sum of allowed embedded stratified simplices.  Notice that in
this situation it makes no sense to talk of $(i+1)$ simplices, for $i$ is
not an integer.

\begin{definition}
Let $i=(i_0,\ldots,i_r)$ be a sequence of non-negative integers.  The
(unrestricted) \emph{local-global homology group} $H_{i;0}X$ is defined to
be the $i$-cycles modulo the $i$-boundaries, as defined in the previous
paragraph. 
\end{definition}

Just as a cone has an apex, so an order~$r$ stratified simplex will have
$r$ apex locii.  These induce a stratification of $\sigma$.  Further
local-global groups can be obtained by imposing conditions on how the
strata of $\sigma$ meet the strata of $X$.  The following definitions make
these concepts clear.

\begin{definition}
Suppose $\sigma=\sigma_i\times C\sigma'$ is an order~$r$ stratified
simplex.  The \emph{$r$-th apex locus} $A_r\sigma$ of $\sigma$ consists
of $\sigma_i\times\{0\}$, where $\{0\}$ is the apex of $C\sigma'$.  The
\emph{$j$-th apex locus} $A_j\sigma$, for $r>j\geq 1$, consists of
$\sigma_i \times C(A_j(\sigma'))$, or in other words $\sigma_i\times
C$ applied to the $j$-th apex locus of $\sigma'$.  We define $A_0\sigma$
to be $\sigma$ itself, but do not count it as an apex locus.  The
$A_j\sigma$ provide a descending filtration of $\sigma$.  The
successive complements $A_j\sigma \setminus A_{j+1}$ are called the
\emph{strata} $A^\circ_j\sigma$ of $\sigma$.
\end{definition}

\begin{definition}
Suppose $f:\sigma\to X$ is an embedded order $r$ stratified simplex.  The
\emph{where} or \emph{$w$-sequence} $w(f)=(w_1,\dots,w_r)$ of $f$ is
defined in the following way.  The entry $w_i$ is defined to be the
largest $j$ such that $f^{-1}(S_j\cup \dots \cup S_n)$ meets the apex
locus $A_i\sigma$ in a set that is dense in $A_i\sigma$.
\end{definition}

\begin{definition}
Let $w=(w_1,\ldots,w_r)$ be a $w$-sequence.  Say that an embedding
$f:\sigma\to X$ of an order $r$ simplex is \emph{$w$-allowed} if for
$1\leq j \leq r$ the inequality $w_j \leq w_j(f)$ holds.
\end{definition}

\begin{definition}
Let $i$ be an index for order $r$ stratified simplices, and let
$w=(w_1,\dots,w_r)$ be a $w$-sequence.  The \emph{$(i,w)$-th local-global
intersection homology group} $H_{i;w}X$ is the quotient group of
$i$-cycles $\eta$ modulo $i$-boundaries $\psi=d\xi$, where each of $\eta$,
$\psi$ and $\xi$ uses only $w$-allowed stratified simplices.  One can also
calls these groups, for $r>0$, the \emph{higher order intersection
homology groups}.
\end{definition}

At the end of \S4 we proved the following result.

\begin{theorem}
If $i$ is of the form $(i_1,0)$, and $w$ is $(w_1)$, then the $(i,w)$
local-global homology $H_{i;w}X$ of $X$ is independent of the
stratification of $X$. 
\end{theorem}

\section{Betti numbers and convex polytopes}

One of the most important properties of middle perversity intersection
homology is that in many cases there is a combinatorial formula for the
Betti numbers.  Here we will say a few words about this, and its relation
to local-global homology.  The basic idea is that from the Betti numbers
for a small class of examples, we can reconstruct the homology theory.

First, we need some algebraic varieties.  If $X$ is an algebraic variety
we will use $IX$ to denote the product of $X$ with $\bfP_1$, while $CX$
will denote the projective cone on $X$.  (This assumes that $X$ is given
to us as a projective variety.)  Now let $W$ be a word of length $n$ in
$I$ and $C$, terminated by a $\{\mbox{pt}\}$.  Thus $W$ denotes an
$n$-dimensional projective variety.  These \emph{$IC$ varieties} are in
fact toric varieties, of rather a special type.

The mpih Betti numbers, also known as the \emph{$h$-vector}, of $IW$ can
be computed from that of $W$ by the K\"unneth formula.  In other words
\begin{equation}
\label{eqn:I}
    h(IW) = IhW
\end{equation}
where the $I$ on the right hand side represents convolution with $(1,1)$,
the Betti numbers of $\bfP_1$.  (Because the odd Betti numbers are here
all zero, it is convenient in this section to omit them.)

Now consider the cone $CW$ on $W$.  It follows easily from the standard
formula for toric variety Betti numbers that $h(CW)$ is obtained from $hW$
by repeating the middle term.  In other words, if $C(hW)$ denotes $hCW$
then
\begin{equation}
\label{eqn:C}
\begin{array}{rl}
C(a)&=(a,a)\\
C(a,a)&=(a,a,a)\\
C(a,b,a)&=(a,b,b,a)\\
C(a,b,b,a)&=(a,b,b,b,a)\\
C(a,b,c,b,a)&=(a,b,c,c,b,a)
\end{array}
\end{equation}
and so on.  It is easy to guess that this should be the required formula,
although proving its truth is another matter.

We can now apply results in combinatorics.  Let $\Delta$ be a convex
polytope.  It has a flag vector $f\Delta$ \cite{bib.MB-LB.gDS,bib.JF.MVIC}.
If $h\Delta$, the Betti numbers of the associated toric variety
$\bfP_\Delta$, is a linear function of $f\Delta$, then it is determined by
its value on the $IC$ polytopes.  This holds for any linear function of
the flag vector, and not just the mpih $h$-vector.  (Here, $I\Delta$ is
the product of $\Delta$ with an interval, a \emph{cylinder} or
\emph{prism}; while $C\Delta$ is the \emph{cone} or \emph{pyramid} on
$\Delta$.  The result of applying a word in $I$ and $C$ to the point
polytope is an $IC$ polytope.) Thus, the formulae (\ref{eqn:I}) and
(\ref{eqn:C}) determine $h\Delta$ for all $\Delta$, provided we either
know or assume that $h\Delta$ is indeed a linear function of $f\Delta$.

Now pretend that we know nothing about intersection homology, except that
there may be a topological homology theory whose Betti numbers are given
by (\ref{eqn:I}) and (\ref{eqn:C}).  With a certain amount of work it is
possible to produce an explicit recursive formula for $h\Delta$ in terms
of $f\Delta$, and when this is done each flag $\delta$ on $\Delta$ makes a
numerical contribution $\lambda_\delta$ to $h\Delta$.  These numerical
contributions can be interpreted in terms of the linear algebra associated
to $\delta$ \cite{bib.JF.CPLA}, and then in terms of the topology of $\PD$.
The conclusion of this process is that, at least conjecturally, the
formulae (\ref{eqn:I}) and (\ref{eqn:C}) give rise to mpih cycles and
boundaries on $\PD$, and hence by way of example tell us what the mpih
conditions on embedded simplices are.

Put in a nutshell, the argument is that (\ref{eqn:I}) and (\ref{eqn:C}),
together with a certain facility in the combinatorics of $\Delta$ and the
topology of $\PD$, would have allowed us to discover middle perversity
intersection homology, if only we had thought to follow this path before
Goresky and MacPherson discovered it via topology.

In this paper we have given a topological definition of local-global
intersection homology.  The actual process of discovery was as described
in this section.  First, rules analogous the the rules (\ref{eqn:I}) and
(\ref{eqn:C}) for $I$ and $C$ were hypothesized, and then the definitions
were unwound to give the topological definition.

The first step was not easy.  Not every pair $I$ and $C$ of rules will
give rise to a linear function on the flag vector.  In fact, they must
satisfy the $IC$-equation \cite{bib.JF.MVIC}
\[
    I\> (I-C)C \> = \> (I-C)C \> I
\]
and also the boundary condition $I\{\mbox{pt}\}= C\{\mbox{pt}\}$, and
these conditions are sufficient.

The reader might wish, as an exercise, to compute perhaps only
heuristically the local-global Betti numbers for $IC$ varieties, perhaps
only in smallish dimensions.  When doing this it is very helpful to know,
as a consequence of a fundamental result of Bayer and Billera
\cite{bib.MB-LB.gDS}, that one wishes for projective toric varieties of
dimension $n$ to have $F_{n+1}$, the $(n+1)$-st Fibonacci number, linearly
independent Betti numbers.  The author's answer to this exercise is the
starting point for \cite{bib.JF.LGIH}, which contains a fuller and less
formal, exposition of concepts presented in this paper.

\section*{Acknowledgements}

The main definition of this paper is the outcome of a long search that was
prompted by the publication in 1985 of \cite{bib.MB-LB.gDS}.  Of the many
people who showed a tolerant interest in the author's previous attempts,
special thanks are due to Marge Bayer, Gil Kalai, Carl Lee, Peter McMullen
and Richard Stanley.  The term `stratified simplex' is due to an anonymous
referee of \cite{bib.JF.CPLA}.

\end{document}